\documentclass[11pt,oneside,reqno]{amsart}
\usepackage{amssymb,epsfig,amsmath,amsthm,amsopn,amstext,amsfonts,amsbsy,graphicx}
\usepackage{color}
\usepackage[latin1]{inputenc}
\usepackage[T1]{fontenc}
\usepackage[english]{babel}
\input epsf

\newtheorem{theorem}{Theorem}
\newtheorem{proposition}{Proposition}
\newtheorem{lemma}{Lemma}

\newtheorem{corollary}{Corollary}

\title[Stochastic Combustion Model]{Fluctuations of the Front in a
  Stochastic Combustion Model \\ $\cdot$ \\
(Fluctuations du Front dans un Mod\`ele de Combustion)
}

\author{Francis Comets$^{1}$, Jeremy Quastel$^{2}$ and Alejandro F. Ram\'\i rez$^{3}$}

\thanks{ AMS 2000 {\it subject classifications}. Primary  82C22, 82B41;
 secondary 82B24, 60K35, 60G99.}


\thanks{$^1$Partially supported by CNRS, UMR 7599.}

\thanks{$^2$Partially supported by NSERC, Canada}

\thanks{$^3$Partially supported by Fondo Nacional de Desarrollo Cient\'\i fico
y Tecnol\'ogico grant 1020686}

\address[Francis Comets]{Laboratoire de Probabilit\'es et Mod\`eles
 Al\'eatoires\\
Universit\'e Paris 7- Denis Diderot\\
2, Place Jussieu\\
F-75 251 Paris Cedex 05, France}

\address[Jeremy Quastel]{Department of Mathematics and Statistics\\
University of Toronto\\
100 St. George Street\\
Toronto, Ontario M5S 3G3, Canada}

\address[Alejandro F. Ram\'\i rez]{Facultad de Matem\'aticas\\
Pontificia Universidad Cat\'olica de Chile\\
Vicu\~na Mackenna 4860, Macul\\
Santiago, Chile}

\bigskip

\email{comets@math.jussieu.fr, quastel@math.toronto.edu, aramirez@mat.puc.cl}

\begin{document}

\begin{abstract}  We consider an interacting particle system on the one dimensional
lattice $\bf Z$ modeling combustion. The process depends on two integer parameters $2\le a<M<\infty$.
Particles move independently as continuous time
simple symmetric random walks except that 1. When a particle jumps to
a site which has not been previously visited by any 
particle, it branches into $a$ particles; 2. When a 
particle jumps to a site with $M$ particles, it is annihilated.
We start from a configuration where all sites to the left of
the origin have been previously visited and study the law of large numbers and
central limit theorem for $r_t$, the rightmost visited site at time $t$.
 The proofs are based on the construction of a
renewal structure leading to a definition of regeneration times for
which good tail estimates can be performed.

\medskip

\noindent
R{\sc \'ESUM\'E.}
On consid\`ere un syst\`eme de particules en interaction sur $\bf Z$
mod\'elisant les particules incandescentes d'un m\'ecanisme de
combustion. Le processus d\'epend de
deux param\`etres entiers $2\le a<M<\infty$. Les particules se
d\'eplacent ind\'ependamment selon des promenades al\'eatoires simples
sym\'etriques \`a temps continu, mises \`a part les interactions suivantes:
1- quand une particule saute vers un site qui n'a jamais encore
\'et\'e visit\'e, elle branche et fait place \`a $a$ particules; 
2- quand une particule saute vers un site abritant 
$M$ particules, elle disparait.
On d\'emarre d'une configuration o\`u seuls les sites \`a gauche de
l'origine ont d\'eja \'et\'e visit\'es et on \'etudie la loi des
grands nombres et le th\'eor\`eme de la limite centrale pour $r_t$,
la position du site le plus \`a droite visit\'e \`a l'instant $t$. 
Les preuves reposent sur la construction d'une structure de
renouvellement associ\'ee \`a des temps de r\'eg\'en\'eration 
dont les queues peuvent \^etre convenablement estim\'ees.

\end{abstract}

\maketitle

\section{Introduction}

The method of regeneration times has been very successfully  applied to problems of random walk in random
environment (for example 
see \cite{s1}, \cite{sz}, \cite{cz}). In a one-dimensional setting  it
was already used via a renewal structure in \cite{kesten}. 
On the other hand, a concept similar to that of regeneration times,
known as {\it cluster structure}, has been well know for many years
in the context of mechanically interacting one-dimensional dynamical systems
of particles
  \cite{sinai} (also see \cite{vladas} and references therein).
 In this article, we extend the method of regeneration times
to an interacting particle system. 
  The system we consider
is a one dimensional stochastic model of combustion.  Heat particles move by symmetric nearest neighbor 
random walks on sites $\{x\in \mathbb Z~:~
x\le r\}$ of the integer lattice, with $r$ representing the position of a flame {\it front}. To the right of the flame front
is a {\it propellant}.  The first heat particle which reaches the propellant at $r+1$ immediately branches into $a\ge 2$
particles, and the front moves one step to the right.  We also include an upper bound $M$ on the number
of particles at each site, so that if a particle tries to jump to a site with $M$ particles, it is immediately killed.
We show that the front in the model moves ballistically to the right, or, more precisely, prove a law of large
numbers for $r_t$.  The next question one might ask is of the fluctuations about the law of large numbers.  We
prove here that these are Gaussian, with a central limit theorem for $t^{-1/2} (r_t - vt)$.

For $M=\infty$ a discrete time version
of the system we are considering has appeared in the literature under the enigmatic name ``frog model'',
and laws of large numbers for the position of the front were proved, also in higher dimensions, using methods
based on the sub-additive ergodic theorem \cite{amp1}, \cite{amp2},
\cite{ampr}, \cite{afm}, \cite{p1}, \cite{p2}.
The continuous time case under the name {\it stochastic combustion
process} was treated in \cite{rs1}, \cite{rs2}.
Note that the $M=\infty$ case is sub-additive, but the $M<\infty$ case is not.  Our main purpose here is to develop
new methods to study such models, and we are especially interested in the fluctuations of the fronts.
If $n_t$ represents the number of particles at the front $r_t$ at time $t$, then $r_t$ moves to the right at rate
$n_t$, and hence $r_t = \int_0^t n_s ds + M_t$ where $M_t$ is a martingale.  The standard approach to the
law of large numbers would then be to show that, as observed from the front, the system has an ergodic invariant
measure $\mu$, and $t^{-1} \int_0^t n_s ds \to E_\mu[n]$.  The central limit theorem would be proved by showing
that the time correlations of $\eta(s)$ decay fast enough.  However, very few methods exist for proving uniqueness,
or ergodicity, of invariant measures of such systems.  So such approaches appear to have limited applicability.
Furthermore, although it seems intuitive that the correlations of $\eta(s)$ decay quickly--likely exponentially fast--
it is not at all apparent how to prove it.   In a future article, we will consider the
case $M=\infty$ which can be handled by essentially the methods developed here, though with major technical
modifications.

The combustion process that we consider is related to deterministic reaction-diffusion equations of the
form $\frac{\partial u}{\partial t} = \frac{\partial^2 u}{\partial x^2}  + f(u)$.  The two key differences are the
discreteness of the variable--we have a number of particles as opposed to a continuum variable $u$--and
the stochasticity.  The effect of discreteness and/or stochasticity on the traveling waves of reaction-diffusion
equations is a question that has not received the attention it deserves, as these effects are likely present in real
systems.  Therefore, we believe it is crucial to develop methods to study such systems.  In the literature
of reaction-diffusion equations one separates several cases according to the behavior of $f$ near zero. 
An $f$ which vanishes on $[0,\theta]$ for some $\theta>0$ with $f(u)>0$ for $u>\theta$ is said to have 
a combustion nonlinearity with ignition temperature cutoff $\theta$ (see, for example \cite{f}).  Note that the discreteness of the particle models of the type we consider makes
them essentially of the combustion type with an ignition temperature of one particle.

In the second section of this article, we  define the combustion process
and state the main result (Theorem \ref{theorem1}). In Section 3, we
define an auxiliary and a labelled process, which will be needed to define the renewal
structure leading to the regeneration times. 
The labelled process can be understood as a
 combustion process where particles are labelled so that if at a given
site a particle has to be killed, it is the one with the smallest
label. The auxiliary process
moves ballistically to the right and is coupled to the labelled
 combustion process
in such a way that it is always to the left of the right-most visited site
$r_t$. In Section  \ref{renewal},
these processes are used to define the renewal structure,
 and the corresponding
regeneration times.
 Then in Section \ref{expectation},
 it is proved that  the regeneration times have
finite   moments of order 2, under appropriate assumption on the
threshold $M$. In Section \ref{limits} we complete the proof of Theorem \ref{theorem1}.

\section{Combustion process}

We define a stochastic process describing the dynamics of particles
on the lattice ${\mathbb Z}$ which move as rate $2$ continuous time simple symmetric
random walks and branch and are killed at  rates depending on the
configuration of neighboring particles.
The branching and killing depend on natural number parameters $2\le a \le M$.
 Each particle performs a continuous time symmetric simple
random walk independent of the others, with two twists:
There is a position $r\in\mathbb Z$ which is the {\it rightmost visited site}.    When a particle jumps right from this site, it branches into $a$ particles, with the
result that there are 
$a$ particles at the new rightmost visited site, $r+1$.
In addition, we allow at most $M$ particles at any site, and maintain this requirement by killing any particle
that attempts to jump to a site with $M$ particles.

The state of our system is

 $$
\Omega: = \{ (r,\eta) :  r \in \mathbb Z,  \eta\in \{ 0,\ldots, 
M\}^{\{ \ldots, r-1,r\} }  \}.
$$
The infinitesimal generator is
\begin{equation}
\nonumber
Lf(r,\eta) = \sum_{x\le r, y \le r , |x-y|=1} \eta(x) 
 (f(r,T_{xy}\eta) -f(r,\eta))\end{equation}
\begin{equation}
\label{by}
 + \eta(r) 
(f(r+1, \eta-\delta_r + a\delta_{r+1}) - f(r,\eta)).
\end{equation}
where $\delta_x$ denotes the configuration with one particle at $x$ and 
\begin{equation}
T_{xy}\eta = \eta-\delta_x + \delta_y 1(\eta(y)<M).
\end{equation}
In the following we will assume  that \begin{equation}
8< M <\infty.
\end{equation}
 We can now state
the main results. Throughout the sequel we will use the notation
$\eta(t,x)$ to denote the number of particles at time $t\ge 0$ at
site $x$ in the stochastic combustion process.

\medskip

\begin{theorem}
\label{theorem1}  Suppose the process is started with initial conditions $r=0$,
$\eta(0,0)\in \{1,\ldots,M\}$ and $\eta(0,x) \in\{0,\ldots,M\}$ arbitrary for $x<0$.

(i) (Law of large numbers).  
There exists  $v$, $0<v<\infty$, which does not
depend on $\eta$, such that a.s.,
\begin{equation}
\lim_{t\to\infty}{r_t}/{t}=v.
\end{equation}
(ii) (Central limit theorem). 
There exists  $\sigma$, $0<\sigma<\infty$, which does not depend on
the initial condition
$\{\eta(0,x):x\in {\mathbb Z}\}$, such that 
\begin{equation}
\label{process}
\epsilon^{1/2}\left(r_{\epsilon^{-1}t}-\epsilon^{-1} vt\right),
\qquad\qquad t\ge 0,
\end{equation}
converges in law as $\epsilon\to 0$  to Brownian motion
with variance  $0<\sigma^2<\infty$.
\end{theorem}
\medskip

\noindent To prove this theorem we will define a renewal structure for the right-most
visited site $r_t$,
in terms of regeneration times, and then will show that such times have finite second
moments.
\medskip

\section{ Auxiliary and labelled processes}

In this section, we will define an auxiliary process
and a labelled process, both coupled to the combustion process.
The auxiliary process $\{\tilde r_t:t\ge 0\}$, will take
values on ${\mathbb Z}$ and will have two fundamental
properties: under certain initial conditions on the combustion
process such that $\tilde r_0=r_0$, 
it will always be to the left of the right-most
visited site $\{r_t:t\ge 0\}$ and
 its dynamics is independent of the particles to the left of
$r_0$ in the combustion process if the appropriate rule
for killing is defined. These properties 
will help us  define the renewal structure 
for the right-most visited site $\{r_t:t\ge 0\}$ of the
combustion process
  and will give us easily a
 lower bound for the limiting speed $v$ of $r_t$.
This is the content of Lemma \ref{coupling}.
The labelled process, corresponds to a process with a state
space larger than the combustion process, where an explicit rule
for killing is given: particles are labelled, and every time a particle
has to be killed, it is the one with the smallest label which
is killed.

Let us first introduce some general notation.
Define ${\mathcal A}:=\{1,2,\ldots ,a\}$,
$B:={\mathbb Z}\times{\mathcal A}$ and
 consider a set of independent continuous time simple symmetric
random walks $\{Y_{x,i}:(x,i)\in B\}$, each one of total jump rate $1$,
and such that $Y_{x,i}(0)=x$ for each $i\in{\mathcal A}$.
 We want to endow the set of indices
$B$ with the lexicographic order.

 Consider a set of independent continuous time rate $2$ simple symmetric
random walks $Y_{x,i}$, with label $(x,i)$,   where $x \in \mathbb Z$ and $i\in \{1,\ldots, a-1\} $. Each random walk $Y_{x,i}$ starts at site $x$
for every $i\in\{1,\ldots , a-1\}$.
 We order the labels by
\begin{equation}\label{lex}
(x,i)  < (x', i') \qquad {\rm if} \qquad x<x' \qquad {\rm or}\qquad x= x' ~~{\rm and}~~ i<i'. 
\end{equation}

\subsection{ Auxiliary process.}
Let $r\in{\mathbb Z}$ and for each $z\in \mathbb Z$, define $B_z$ as the set of
 $M$ labels not exceeding $(z,a-1)$ in the order (\ref{lex}).  Define
$A_{r,z}$ to be the labels in $B_z$ with $x\ge r$.

We define a sequence of waiting times.
Let $\nu_0:=0$ and  define $\nu_1$ as the first time  one of the 
random walks $\{Y_{z,i}: (z,i)\in A_{r,r}\}$,
hits the site $r+1$.  Next, define
$\nu_2$ as the first time one of the random walks $\{Y_{z,i}:(z,i)\in
A_{r,r+1}\}$,
hits the site $r+2$.  In general,   for $k\ge 3$, we define
 $\nu_{k}$ as the first
 time  one of the random walks $\{Y_{z,i}:(z,i)\in A_{r,r+k-1}\}$,
visits the site $r+k$.
Finally,  for $n\in{\mathbb N}$, let
\begin{eqnarray}
\label{let}
&\tilde r^r_t:=r+n,\qquad {\rm if}\qquad \sum_{k=0}^n\nu_k\le t<\sum_{k=0}^{n+1}\nu_k,&
\end{eqnarray}
where the superscript $r$ in $\tilde r^r_t$ indicates that $\tilde r^r_0=r$.

\begin{lemma}
\label{lln-auxiliary}
There exists an $\alpha>0$ which does not depend on $r$
such that a.s.
\begin{equation}\label{alpha}
\lim_{t\to\infty}{\tilde r^r_t}/{t}=\alpha.
\end{equation}
\end{lemma}

\begin{proof}[Proof] Once one notes that for every $j\ge 0$ the 
random variables
 $\{\nu_{k\ell +j}
:k\ge 1\}$ are independent whenever $\ell>M/a$, the proof is a simple exercise.
\end{proof}
\medskip

\subsection{ Labeled process.} We enlarge the state space of the stochastic 
combution process so that particles carry labels which tell us where they
originated.

We will want to keep track of where particles came from, even after restarting at stopping times.
 Hence each particle will have a 
starting position $z\in\mathbb Z$ and label $(x,i)$, $x\in\mathbb Z$, $i\in \{1,\ldots,a-1\}$ describing
its birthplace.  We will allow the possibility that $z\ne x$.

  At time $0$, we have an $r\in\mathbb Z$
representing the rightmost visited site, and a subset $\mathcal I(0)$ of the labels $(x,i)$ with $x\le r$, representing the set of labels of live particles
at time $0$.  
To each one of these labels is assigned a position $z=Z_{x,i}(0)\le r$ which is the position at time $t=0$ of
that particle.  The position at time $t$ is $Z_{x,i}(t) = Y_{x,i} (t) + z- x$.

 To keep track of the killing in this process, we make a rule that 
whenever a particle jumps to a site with $M$ particles,
the particle at that site with the smallest label is killed and
the corresponding label is removed from the set $\mathcal I$ of labels
of live particles.

When a particle jumps to site $r+1$, the labels $\{(r+1,1),\ldots,
(r+1,a-1)\}$  with corresponding to particles
with initial positions $r+1$, are added to the set of labels of live particles.
The time this happens will be denoted
 $\rho_1$.
  These particles then have trajectories $Z_{r+1,i}(t)$ which is equal to $Y_{r+1,i }(t-\rho_1)$ for $t\ge\rho_1$.

Similarly, for $k\ge 2$,
 $\rho_1+\cdots+\rho_k$ will be the first time a particle jumps to $r+k$ and at that time $\{(r+k,1),\ldots,
(r+k,a-1)\}$ are added to $\mathcal I $,
 with trajectories $Z_{r+k,i}(t) = Y_{r+k,i}(t-\rho_1-\cdots-\rho_k)$ for 
$t\ge \rho_1+\cdots+\rho_k$ until such time as the particles are killed and their label removed from the set of labels of live particles.

We denote by $\mathcal I(t)$ the set of labels of live particles at time $t$ and by ${\mathcal Z}(t)=\{Z_{(x,i)}(t): (x,i)\in{\mathcal I}(t)\}$
the positions of the corresponding random walks.   

To avoid pathologies it is useful to insist that initially the 
set of labels of live particles includes at least one with $x=r$.
 The rightmost visited
site $r$ is the supremum of the $x$ over the collection of labels,
$r_t = \sup\{ x~:~ (x,i) \in \mathcal I(t)\}$. 
Let us formalize the above discussion. Call ${\mathbb L}$ the triples $(r,{\mathcal I},
{\mathcal Z})$ formed by an integer $r\in {\mathbb Z}$, a set of
labels ${\mathcal I}\subset\{(x,i):x\le r, 1\le i\le a-1\}$ and
position function ${\mathcal Z}:{\mathcal I}\to \{\ldots, r-2,r-1,r\}$ taking
values in the integers smaller than or equal to $r$. We now define

\begin{equation}
\nonumber
{\mathbb S}:=\left\{(r,{\mathcal I},{\mathcal Z})\subset{\mathbb L}:
\max_{(x,i)\in{\mathcal I}}x=r, \max_{z\le r}\sum_{(x,i)\in{\mathcal I}}{\bf 1}_{{\mathcal Z}(x,i)}(z)\le M
\right\}.
\end{equation} 
which will be the state space of our process.

The {\it labeled process} starting from $w=(r,{\mathcal I}(0),{\mathcal Z}(0))$,
is now defined as the triple $w_t=\{(r_t,{\mathcal I}(t),{\mathcal Z}(t)):t\ge 0\}$,
defining a strong Markov process taking values on ${\mathbb S}$, and with a law
given by a probability measure ${\mathbb P}_w$ defined on the Skorohod space $D([0,\infty);{\mathbb S})$.

The {\it combustion process} described in section 1 is the particle count
$\eta(t):=\{\eta(y,t):y\le r_t\}$ where
$\eta(y,t) = \sum_{(x,i)\in{\mathcal I}(t)} 1(Z_{x,i} (t)=y)$.
 We will 
occasionally use the more explicit notation
notation $\eta_w(t)$ and $\eta_w(y,t)$ instead of $\eta(t)$
and $\eta(y,t)$ respectively, with the understanding that
$\eta_w(0)$ is the particle count of the initial condition $w$.

We already defined $\rho_1$ to be the first time that one of these particles hits $r+1$.

\begin{lemma}
  Suppose that  $(r,1),\ldots, (r,a-1)\in \mathcal I(0)$, all initially at $r$.  Then $\rho_1\le \nu_1$.
\end{lemma}

\begin{proof}
The dynamics of the $a$ particles which start at $r$ is the same in the labeled process as in the ones we
look at in the auxiliary
process.  Indeed, these are the only particles considered in the auxiliary process up to time $\nu_1$, while the labeled
process may have many others, each of which has a chance to be the first to hit $r+1$.
\end{proof}

Recall that
 $\rho_k+\cdots+\rho_1$ is the first time one of the labeled particles
in the labeled process hits $r+k$.
  Note of course that the hitting could be done
by one of the $a$ new particles created at $r+i$ at time $\rho_i$, $i<k$.

\begin{lemma}
\label{difficult}
Suppose that $(r,1),\ldots, (r,a-1)\in \mathcal I(0)$, all initially at $r$.
  Then $\rho_k\le \nu_k$.
\end{lemma}

\begin{proof}  Before giving a general proof, we describe the special case of $M=a=k=2$ where the idea
is more transparent.

Case $M=a=k=2$:  Note first of all that if $a=2$, we have $i=1$ always and hence we can drop the $i$ in the
labels.  The labeled process starts with one particle with label $r$ at $r$, possibly one other particle at site
$r$, with a label smaller than $r$, and other particles with labels 
smaller than $r$ at arbitrary positions to the left of $r$.
At time $\rho_1$ we have one particle labeled $r+1$ at $r+1$, which, after that time, has trajectory $Z_{r+1}(t)=
Y_{r+1}(t-\rho_1)$,
and another particle, labeled $r$ at some site $x\le r+1$, which, after that time, has trajectory $Z_r(t)= Y_{r}(t)$.
Note that neither particle can be killed before time $\rho_1+\rho_2$ because they are the two with the highest 
labels until that time.  There could also be other particles with other labels at positions $x\le r$.  Denote by
$\tau_{r+1}=\inf\{ t\ge \rho_1 : Z_{r+1}(t)=r+2\}$, by $\tau_r=\inf\{t\ge \rho_1: Z_r(t)=r+2\}$ and
by $\tau_{others}$ the first time one of the others hits.  Then $\rho_2=  \min\{ \tau_{r+1},\tau_r,\tau_{others}\}-\rho_1
\le \min\{ \tau_{r+1},\tau_r\}-\rho_1$.

On the other hand, $\nu_2= \min\{ \sigma_r, \sigma_{r+1} \}$ where $\sigma_r= \inf\{ t\ge 0: Y_{r}(t)=r+2\}$
and $\sigma_{r+1}=\inf\{ t\ge 0: Y_{r+1}(t) = r+2\}$.  Now
\begin{eqnarray} &
\tau_{r+1} = \inf\{ t\ge \rho_1 : Z_{r+1}(t) = r+2\}  = \inf\{ t\ge \rho_1 : Y_{r+1}(t-\rho_1) = r+2\} &\nonumber
\\
&
 = \inf\{ t\ge 0  : Y_{r+1}(t) = r+2\}+\rho_1 = \sigma_{r+1}+\rho_1. &\end{eqnarray}
and
\begin{eqnarray} 
\nonumber &
\tau_{r} = \inf\{ t\ge \rho_1 : Z_{r}(t) = r+2\}  = \inf\{ t\ge \rho_1 : Y_{r}(t) = r+2\} &
\\
&
 \le  \inf\{ t\ge 0  : Y_{r}(t) = r+2\}+\rho_1 = \sigma_{r}+\rho_1. \nonumber&\end{eqnarray}
The result follows.

General case:
Let us examine the configuration in the labeled process at time $\rho_1+\cdots+\rho_{k-1}$.  
It has $a-1$ particles at site $r+k-1$  with labels $(r+k-1,1), \ldots,(r+k-1,a-1)$ plus one additional
particle with an unknown label, the one which hit $r+k-1$.  An additional $M-a+1$ particles (which might include
the additional particle of unknown label which hit $r+k-1$) with the
next highest labels after the first $a-1$ have some unknown positions, and there may in addition be any number of other particles in the configuration..  In the time interval $[0,\rho_1+\cdots+\rho_{k-1}]$
none of the $M-a+1$ particles has hit $r+k$.   Their trajectories in the time interval 
$[\rho_1+\cdots+\rho_{k-1},\rho_1+\cdots+\rho_{k}]$ are $Y_{(r+k-2,a-1)}(t-\rho_1+\cdots+\rho_{k-2}),\ldots$.
The definition of $\nu_k$, on the other hand, involves the particles with these $M$ leading indices, but without the
time shift.  Hence the first time that any of the first $a-1$ particles hits is identical to that of the auxiliary process,
the first time that any of the next $M-a+1$ particles hits is greater in the auxiliary process, and one of the other
possible particles in the labeled process could be the one which hit first, making the time shorter still.
\end{proof}

\subsection{ Enlarged process.}  We have seen that the labeled process is defined in terms of a set of
labels ${\mathcal I}(t)$, for $t\ge 0$, with the property $r_t=\sup\{x: (x,i)\in{\mathcal I}(t)\}<\infty$. Whenever
a particle with a label in ${\mathcal I}(t)$ is killed, this label is removed.
In the enlarged process we 
keep track of the killed particles as well.
 Define for each $t\ge 0$,
the set of labels of all activated particles up to time $t$ in the labeled process,
\begin{equation}
\nonumber
\bar{\mathcal I}(t)=\cup_{0\le s\le t}{\mathcal I}(s).
\end{equation}
Consider the set $\bar{\mathcal Z}(t)=\{Z_{x,i}(t):(x,i)\in \bar{\mathcal I}(t)\}$ of all the corresponding
random walks at time $t$. The {\it enlarged process} 
is defined as the triple
 $\bar w_t=\{(r_t,\bar{\mathcal I}(t),\bar{\mathcal Z}(t)):t\ge 0\}$. 

Consider now the set of labels ${\mathcal R}(t)$, obtained after removing from ${\mathcal I}(t)$ all
labels $(x,i)$ with $x<r=\sup\{y:((y,i)\in {\mathcal I}(0)\}$. 
We define for $y\le r_t$ the particle count \begin{equation}
\zeta(y,t)=\sum_{(x,i)\in{\mathcal R}(t)} 1(Z_{(x,i)}(t)=y).\end{equation}

Also let ${\mathcal L}(t)$ be the set of labels obtained after removing from ${\mathcal I}(t)$ all
labels $(x,i)$ with positions $x\ge r$. 
We define for $y\le r_t$ the particle count \begin{equation}
\phi(y,t)=\sum_{(x,i)\in{\mathcal L}(t)} 1(Z_{(x,i)}(t)=y).\end{equation}

We similarly define $\bar{\mathcal L}(t)$ as the set of labels obtained after
removing from $\bar{\mathcal I}(t)$ all labels $(x,i)$ with positions $x\ge r$, and the corresponding
particle count for $y\le r_t$ as 
\begin{equation}\bar\phi(y,t)= \sum_{(x,i)\in\bar{\mathcal L}(t)}1(Z_{x,i}(t)=y).\end{equation}

\begin{lemma}
\label{coupling}
 (i) For every initial condition $w\in {\mathbb S}$ and every $t\ge 0$ and $x\le r_t$
\begin{equation}\phi(x,t)\le \bar\phi(x,t).\end{equation}

\noindent (ii) For every $w=(r,{\mathcal I}(0),{\mathcal Z}(0))\in{\mathbb S}$ with labels $(r,1),\ldots ,(r,a-1)$
 at site $r$ and corresponding initial positions $Z_{(r,i)}(0)=r$ for $1\le i\le a-1$,
\begin{equation}
\tilde r^r_t\le r_t.
\end{equation}

\noindent (iii) For every $w=(r,{\mathcal I}(0),{\mathcal Z}(0))\in{\mathbb S}$,
 the processes $\bar\phi$
and $\tilde r^r_t$ are independent.
\end{lemma}

\begin{proof} Part $(i)$ follows directly from the definitions. Part $(ii)$ is a consequence of Lemma \ref{difficult}.
Part $(iii)$ follows from the observation that $\bar\phi$ and $r^r_t$ are defined in terms of the random walks
$\{Y_{x,i}:x<r\}$ and $\{Y_{x,i}:x\ge r\}$ respectively, which are independent.
\end{proof}

\section{The renewal structure}  
\label{renewal}
Consider the enlarged process $\bar w_t$ with its natural
filtration $\bar{\mathcal F}_t$
with an  initial condition $w_0$ having particles with labels $(r,1),\ldots,(r,a-1)$  at site $r$, and any allowable
configuration of particles with labels to the left of $r$.   Let $\alpha=\lim_{t\to\infty} \tilde r_t/t$ be as in (\ref{alpha})
and choose any
$0<\alpha' < \alpha$.  Define
\begin{equation}
\nonumber
U:=\inf \{t\ge  0: \tilde r_t-r <  \lfloor\alpha' t\rfloor \}
\end{equation}
where $\lfloor x\rfloor$ denotes the greatest integer less than or equal to $x$.

Define $V$ as the first time one of the 
particles initially to the left
of $r-1$  hits $\{ x~:~ x\ge \lfloor\alpha' t\rfloor+r\}$,
$$
V:=\inf\{t\ge 0:\sup_{x\ge [\alpha' t]+r}\bar\phi_{w_0}
(x,t)>0\}.
$$  
$U$ and $V$ are stopping times with respect to 
$\{\bar {\mathcal F}_t:t\ge 0\}$.  Let
\begin{equation}
D:=\min\{U,V\}.
\end{equation}
For each $y\in{\mathbb Z}$, let 
$$T_y:=\inf\{t\ge 0:r_t\ge y\}$$
 denote the 
hitting time of $y$ by the rightmost visited site in the labeled  process $w_t$. 

We will also need the first time $U$ and $V$ happen
 after time $s\ge 0$,
\begin{equation}
U\circ\theta_s:=\inf \{t\ge  0: \tilde r^{r_s}_t- r_s <  \lfloor\alpha' t\rfloor \},
\end{equation}
\begin{equation} \label{def:V}
V\circ\theta_s:=\inf\{t\ge 0:\sup_{x\ge [\alpha' t]+r_s}\bar\phi_{w_s}
(x,t)>0\},\end{equation}
and $D\circ\theta_s:=\min\{U\circ\theta_s,V\circ\theta_s\}$.

Choose an integer $L:=M$ and define sequences $\{S_k:k\ge 0\}$ and $\{D_k:k\ge 1\}$
of $\bar{\mathcal F}_t$-stopping times as follows.  Start with 
$S_0=0$ and $R_0=r$.  Then define
\begin{equation}
\nonumber
S_1:=T_{R_0+L}\qquad D_1:=D\circ \theta_{S_1}+S_1,\qquad R_1:=r_{D_1},
\end{equation}
and, for $k\ge 1$,
\begin{equation}
\nonumber
S_{k+1}:=T_{R_k+L}\qquad D_{k+1}:=D\circ\theta_{S_{k+1}}+S_{k+1},\qquad R_{k+1}:=r_{D_{k+1}}.
\end{equation}
Note that these times are not necessarily finite and we make the
convention  that $r_\infty=\infty$.  
Similarly, we define $ U_{k}:=U\circ\theta_{S_{k}}+S_{k}$ and 
$ V_{k}:=V\circ\theta_{S_{k}}+S_{k}$ for $k \geq 1$.
Let
\begin{equation}
\nonumber
K:=\inf\{k\ge 1:S_k<\infty ,D_k=\infty\},
\end{equation}
and define the {\it regeneration time}, 
\begin{equation}
\label{six}
\kappa:=S_{K}.
\end{equation}
Note that $\kappa$ is {\it not} a stopping time with respect to $\bar{\mathcal F}_t$.

Denote by $\mathcal G$ the information up to time $\kappa$, defined as
 the completion with respect to ${\mathbb P}_{w}$ 
of the smallest $\sigma$-algebra containing all sets of the form $\{\kappa\le t\}\cap A$, $A\in\bar{\mathcal F}_t$.

In section \ref{expectation}, we will show that ${\mathbb E}_w[ \kappa^2 ] <\infty$ and hence, in particular
\begin{equation}\label{null}
{\mathbb P}_w[ \kappa <\infty ] =1.
\end{equation}

\medskip
\begin{proposition} 
\label{proposition1} Let $A$ be a Borel subset of 
$D([0,\infty);\Omega)$. Then,
\begin{equation}
\nonumber
 {\mathbb P}_{w}\left[\tau_{-r_\kappa}\zeta(\kappa+\cdot)\in
A ~|~{\mathcal G}\right]
={\mathbb P}_{a\delta_0}[   \eta(\cdot)\in A ~|~ U=\infty],
\end{equation}
where $a\delta_0$ denotes a configuration with $a$ particles at $0$ and none anywhere else.
\end{proposition}

\begin{proof}[Proof] We have to show that for any
 $B\in{\mathcal G}$,
\begin{equation}\label{cc}
{\mathbb P}_{w}[B, \{\tau_{-r_{\kappa}} \zeta(\kappa+\cdot)\in A\} ]
={\mathbb P}_{w}[B]~{\mathbb P}_{a\delta_0}[ \eta(\cdot)\in A ~|~ U=\infty].
\end{equation}
Now, using (\ref{null}),
\begin{eqnarray}
\nonumber
&\ &{\mathbb P}_{w}[B, 
\{ \tau_{-r_{\kappa}} \zeta(\kappa+\cdot)\in A\} ]
={\mathbb P}_{w}[\{\kappa<\infty\}, B, \{ \tau_{-r_{\kappa}} 
\zeta(\kappa+\cdot)\in A\} ]
\\
\nonumber
&=&\sum_{k=1}^\infty {\mathbb P}_{w}\left[\{S_k<\infty, D_k=\infty\}, B,  
 \tau_{-r_{\kappa}}\zeta(\kappa+\cdot)\in A \right]
\\
\label{one}
&=&
\sum_{k=1}^\infty\sum_{x\in {\mathbb Z}}{\mathbb P}_{w}[r_{S_k}=x,
S_k<\infty, D_k=\infty , B,
  \tau_{-x} \zeta(S_k+\cdot)\in A ].
\end{eqnarray}
 From the definition of $\mathcal G$ there is an event $B_k\in\bar{\mathcal F}_{S_k}$
such that
$B_k=B$ on $\kappa=S_k$.
Therefore, we can continue developing (\ref{one}) to obtain,
\begin{eqnarray}
\nonumber
&=&\sum_{k=0}^\infty\sum_{x\in {\mathbb Z}}
{\mathbb P}_{w}\left[r_{S_k}=x, S_k<\infty, D_k=\infty, B_k, \tau_{-x}
\zeta(S_k+\cdot)\in A \right]\\
\nonumber
&=&\sum_{k,x}
{\mathbb E}_{w}\left[{\bf 1}_{r_{S_k}=x, S_k<\infty, B_k}
{\mathbb P}_{w} \left[ D_k=\infty, \tau_{-x}\zeta(S_k+\cdot)\in A
~|~\bar{\mathcal F}_{S_k}\right] 
\right],
\end{eqnarray}
where ${\mathbb E}_w$ is the expectation defined by ${\mathbb P}_w$.
But on the events $S_k<\infty$ and $r_{S_k}=x$, we have by parts $(i)$ and $(ii)$
of Lemma \ref{coupling}, that
\begin{equation}
\label{zeta-eta}
\zeta_{w}(S_k+\cdot )=\eta_{a\delta_{x}}(\cdot)
\end{equation}
when $U_k=V_k=\infty$, and that $\eta_{a\delta_{S_k}}(\cdot)$ is independent of the configuration
of particles to the left of $x$. 
Indeed, part $(i)$ of Lemma \ref{coupling}, and the event $V_k=\infty$,
imply that the particles with initial positions $z$ to the left of $x$,
 never hit the line $\lfloor \alpha' t\rfloor +x$. And
part $(ii)$ of Lemma \ref{coupling}, and the event $U_k=\infty$, imply that
the front $r_t$ is always to the right of the same line. Hence, there is no
effect of the particles initially to the left of $x$ on the front $r_t$,
so that $\zeta_w(S_k+\cdot)=\eta_{a\delta_{x}}(\cdot)$.
Then, (\ref{zeta-eta}) combined with
 the independence of $U_k$ and $V_k$ given $\bar{\mathcal F}_{S_k}$ by part $(iii)$ of
Lemma \ref{coupling},
the translation invariance,  and the strong Markov property imply that on the
events $S_k<\infty$ and $r_{S_k}=x$,
\begin{eqnarray}
\nonumber
 &&
{\mathbb P}_{w} \left[  U_k=\infty,V_k=\infty,\tau_{-x}
\zeta(S_k+\cdot)\in A
~|~\bar{\mathcal F}_{S_k}\right]\\
\nonumber
&=&{\mathbb P}_{w} \left[  U_k=\infty,\tau_{-x}
\eta_{a\delta_{x}}(\cdot)\in A ~|~ \bar{\mathcal F}_{S_k}\right]
 {\mathbb P}_{w}\left[V_k=\infty ~|~ \bar{\mathcal F}_{S_k}\right]\\
\nonumber
& =&
{\mathbb P}_{a\delta_0}[ U=\infty,\eta(\cdot)\in A]{\mathbb P}_{w}[  V_k=\infty ~|~\bar{\mathcal F}_{S_k}] . 
\end{eqnarray}
Summarizing, we have,
\begin{eqnarray}
\nonumber
&&
{\mathbb P}_{w}[\kappa<\infty]~{\mathbb P}_{w}[B,\tau_{-r_{\kappa}} 
\zeta(\kappa+\cdot)\in A ]\\
\label{aa}
&&
=  {\mathbb P}_{a\delta_0}[ U=\infty,\eta(\cdot)\in A]
\sum_{k,x} 
{\mathbb P}_{w}[ V_k=\infty,r_{S_k}=x,S_k<\infty, B_k].
\end{eqnarray}
Letting $A=\Omega$ gives 
\begin{equation}
\label{bb}
{\mathbb P}_{w}[\kappa<\infty] {\mathbb P}_{w}[B]=
{\mathbb  P}_{a\delta_0}[U=\infty] \sum_{k,x} 
{\mathbb P}_{w}[ V_k=\infty,r_{S_k}=x,S_k<\infty, B_k].
\end{equation}
(\ref{aa}) and (\ref{bb}) together imply (\ref{cc}).\end{proof}

  Now define  $\kappa_1\le\kappa_2\le\cdots $ by $\kappa_1:=\kappa$ and for $n\ge 1$
\begin{equation}
\kappa_{n+1}:=\kappa_n + \kappa (\bar w_{\kappa_n+\cdot}).
\end{equation}
where $\kappa(\bar w_{\kappa_n+\cdot})$ is the
regeneration time starting from $\bar w_{\kappa_n+\cdot}$ and
we set $\kappa_{n+1}=\infty$ on $\kappa_n=\infty$ for $n\ge 1$. We will call
$\kappa_1$ the {\it first regeneration time} and $\kappa_n$ the {\it $n$-th regeneration time}.




For each $n\ge 1$ the 
$\sigma$-algebra,
${\mathcal G}_n$
will be the completion with respect to ${\mathbb P}_{w}$ of the smallest $\sigma$-algebra containing all
 sets of the form $\{\kappa_1\le t_1\}\cap
\cdots\cap\{\kappa_n\le t_n\}\cap A$, $A\in{\mathcal F}_{t_n}$.
Clearly ${\mathcal G}_1={\mathcal G}$.

\begin{lemma}
\label{lemma3} $\{U=\infty\}\in{\mathcal G}_1$.
\end{lemma}

\begin{proof}[Proof] Note that $\{\kappa_1=\infty\}$ is a null event for $\mathbb P_w$ and hence, since $\mathcal G_1$
is complete, it is enough to show
that $\{ U<\infty\} \cap \{ \kappa_1<\infty\} \in {\mathcal G}_1$.  

For notational convenience, we write $\tilde r^k_{\cdot}$ instead of $\tilde r^{r_{S_k}}_{\cdot}$.
Note that whenever $U<\infty$, $S_k<\infty$ and
$\tilde r_U> r_{S_k}$ happen, for some
$k\ge 1$, then necessarily  $\kappa_1>S_k$. Indeed, from
the observation that $\tilde r_{S_k+\cdot}=\tilde r_{\cdot}^{k}$,
we see that if the events $U<\infty$, $S_k<\infty$ and $\tilde r_U\ge r_{S_k}$ happen,
 then 
we must have that $U_k<\infty$, and hence $D_k<\infty$.
By summing over the intersection with $\{\kappa_1=S_k\}$ we see that it follows 
 that  $\{U<\infty\}\cap\{\kappa_1<\infty\}\cap\{\tilde r_U> r_{\kappa_1}\}=\emptyset$.
 So
\begin{equation}
 \{U<\infty\}\cap\{\kappa_1<\infty\}=\{\tilde 
r_U\le r_{\kappa_1}\}\cap\{\kappa_1<\infty\}.
\end{equation}
Since  $\tilde r_U\le r_t$ implies that $U<\infty$,
it follows that $\{ \tilde r_U > r_{\kappa_1}\} \cap \{\tilde 
r_U\le r_t\}\cap\{\kappa_1\le t\} $ is empty.  Hence 
$\{\tilde 
r_U\le r_t\}\cap\{\kappa_1<t\} = \{ \tilde r_U\le r_{\kappa_1}\} \cap \{\tilde 
r_U\le r_t\}\cap\{\kappa_1<t\}  = \{ \tilde r_U\le r_{\kappa_1}\} \cap\{\kappa_1<t\}$.   Thus
\begin{equation}
\{\tilde 
r_U\le r_{\kappa_1}\}\cap\{\kappa_1<\infty\}=
\bigcup_{n=1}^\infty \{\tilde 
r_U\le r_n\}\cap\{\kappa_1<n\}.
\end{equation}
The result then follows from the fact that
$\{\tilde  r_U\le r_t\}\in \bar{\mathcal F}_t$ for each $t>0$
which is a direct consequence of the construction of the processes.
\end{proof}

\begin{proposition}
\label{cmp}
  Let $A$ be a Borel subset of 
$D([0,\infty);\Omega)$. Then,
\begin{equation}
\nonumber
{\mathbb P}_{w}[\tau_{-r_{\kappa_n}}\zeta(\kappa_n+\cdot)
\in A  ~|~{\mathcal G}_n]
={\mathbb P}_{a\delta_0}[  \eta(\cdot)\in A ~|~ U=\infty].
\end{equation}
\end{proposition}

\begin{proof}[Proof] 
 Let $\psi:D[[0,\infty);{\mathbb S})
\to D[[0,\infty);{\mathbb S})$ be the map given by,
\begin{equation}
\psi (\bar w)(\cdot )=\tau_{-r_{\kappa_1}} \bar w(\kappa_1 (\bar w(\cdot ))+\cdot ).
\end{equation}
Then note that, ${\mathcal G}_{k+1}$ is generated by ${\mathcal G}_1$
and $\psi^{-1}\left({\mathcal G}_k\right)$, and
that the $\sigma$-algebras ${\mathcal G}_1$ and $\psi^{-1}({\mathcal G}_k)$ are
independent.
 The proof of this theorem now follows from the above observations,
 induction on $n\in{\mathbb Z}^+$ using
Proposition \ref{proposition1}, and Lemma \ref{lemma3}.
\end{proof}

\begin{proposition}
\label{iid}
 Let $w\in {\mathbb S}$.
(i)  Under ${\mathbb P}_w$, $\kappa_1,\kappa_2-\kappa_1, \kappa_3-\kappa_2, \ldots$ are independent, and $\kappa_2-\kappa_1, \kappa_3-\kappa_2, \ldots$ are identically distributed with law identical to that of $\kappa_1$ under
${\mathbb P}_{a\delta_0}[\cdot |U=\infty]$.
(ii) Under ${\mathbb P}_w$, $ r_{\cdot\land\kappa_1}, r_{(\kappa_1+\cdot)\land\kappa_2} - r_{\kappa_1}, r_{(\kappa_2+\cdot)\land\kappa_3} - r_{\kappa_2}, \ldots$ are independent, and $r_{(\kappa_1+\cdot)\land\kappa_2} - r_{\kappa_1}, r_{(\kappa_2+\cdot)\land\kappa_3} - r_{\kappa_2}, \ldots$ are identically distributed with law identical to that of $r_{\kappa_1}$ under
${\mathbb P}_{a\delta_0}[\cdot|U=\infty]$ .
\end{proposition}

\begin{proof}[Proof] This follows directly from Proposition \ref{cmp}.
\end{proof}

\section{Expectations and variances of the regeneration times}
\label{expectation}

\subsection{Estimates for the auxiliary process.} 
In this subsection we obtain some estimates for the auxiliary process.  The process will always start from $a\delta_0$
and we denote the corresponding measure on the trajectories of the auxiliary process by $P$, and expectations by $E$.

\begin{lemma}
\label{moment1}
For $1\le p < M/2$ and $j\ge M$, there exists a constant $C=C(p)<\infty$ such that
\begin{equation} E [ |\nu_j|^p] < C.\end{equation}
\end{lemma}

\begin{proof}  For $j\ge M$ the auxiliary process has its full complement of $M$ particles attempting
to hit $j+1$.
$\nu_j$ is then the minimum of $\gamma_1,\ldots,\gamma_M$ which are the hitting times of $1$ of 
$M$ random walks starting in $\{-M/a \le x\le 0\}$.  Standard estimates for random  walks
give $P[ \gamma_i > t] \le Ct^{-1/2}$ for some $C<\infty$.  Hence $P[\nu_j>t] \le C^M t^{-M/2}$.  
\end{proof}

\begin{lemma} 
\label{est-U}
For $1\le p<M/2$, there exists a constant $C=C(p)<\infty$
such that for all initial conditions $w$, and all $t>0$
\begin{equation}
\label{such-that}
  P\left[ t<U<\infty\right]\le C{t^{-(p/2)}}.
\end{equation}
\end{lemma}

\begin{proof} By translation invariance we can assume $r=0$. Let $t_1$ be such that $\lfloor \alpha't_1\rfloor=M$.
Then, when $t\ge t_1$, we have that 
\begin{equation}
 P\left[ t<U<\infty\right]\le  P \left[\tilde r^0_{t_1}> M,
\cup_{s> t}\left\{\tilde r^0_s\le \lfloor \alpha' s\rfloor\right\}\right].
\end{equation}
 But, if $\tilde r^0_s\le \lfloor\alpha' s\rfloor$
for $s\ge t$,
then $\sum_{j=1}^{\lfloor\alpha' s\rfloor}\nu_j\ge s$,
which in turn implies that,
$\frac{1}{n}\sum_{j=1}^n\nu_j \ge \frac{1}{\alpha'}$ for some 
$n\ge \lfloor \alpha' t\rfloor$. Similarly,
$\tilde r^0_{t_1}>M$, implies that $\sum_{j=1}^M\nu_j\le (M+1)/\alpha'$.
Therefore, whenever $\lfloor \alpha't\rfloor\ge M+1$, we have
\begin{eqnarray}
 P\left[ t<U<\infty\right]\le &
 P \left[\sum_{j=1}^M\nu_j\le \frac{M+1}{\alpha'}, \bigcup_{n=\lfloor \alpha' t\rfloor}^\infty
\left\{\frac{1}{n}\sum_{j=1}^n\nu_j\ge \frac{1}{\alpha'}   \right\}
\right] &\nonumber \\
\le &
 P\left[\bigcup_{n=\lfloor \alpha' t\rfloor}^\infty\left\{\frac{1}{n}\sum_{j=M+1}^n\nu_j\ge \frac{1}{\alpha'}
\left(1-\frac{M+1}{n}\right)\right\}\right]. &
\end{eqnarray}
Now remark that for $j\ge M$, the random 
variables $\nu_j$ are identically distributed and have finite moments
of order $p<M/2$  by Lemma \ref{moment1}. Each has expected value
 $1/\alpha$.  Define $\gamma_j=\nu_j-1/\alpha$.
Choose $t_2$ as any real number such that $\beta=
\frac{1}{\alpha'}-\frac{1}{\alpha}-\frac{M+1}{\lfloor \alpha' t_2\rfloor}>0$. Then, whenever $t\ge t_2$, if  $N= \lfloor \alpha' t\rfloor$,
we have that,
\begin{eqnarray} \label{eq:toto}
P\left[ t<U<\infty\right]\le & 
P\left[\sup_{n\ge N } \frac{1}{n}\sum_{j=M+1}^n\gamma_j\ge \beta\right].&
\end{eqnarray}
Recall that for each $0\le i< \ell$,  $\ell = \lfloor M/a\rfloor +1$, the random variables $\{\nu_{k\ell +i}~:~ k\ge 1\}$ are independent.
Observe from (\ref{eq:toto}) that,  if $\beta'= \beta a/M$,
\begin{eqnarray}
P\left[ t<U<\infty\right]\le & \sum_{i=0}^{\ell-1} 
P\left[\sup_{n\ge N } \frac{1}{n}\sum_{j=M+1, j= k\ell +i }^n\gamma_j\ge \beta'\right]&\label{you}
\end{eqnarray}

Suppose $X_1,X_2,\ldots$ are independent and identically distributed random variables with mean $0$.  By
 Kolmogorov's inequality,
\begin{eqnarray} &
P[ \sup_{n\ge N}\frac{1}{n}\sum_{i=1}^n X_i \ge \epsilon ] \le   
\sum_{k=0}^\infty  P \left[\sup_{ 2^kN\le n\le  2^{k+1}N}
\frac{1}{ 2^kN}\sum_{i=1}^nX_i\ge\epsilon\right] &\nonumber
\\ &
\label{thirtysix}
\le \sum_{k=0}^\infty (2^kN\epsilon)^{-p}  E[|\sum_{i=1}^{2^{k+1}N} X_i|^p] &
\end{eqnarray}
Now, for $p\ge 2$, if $E[|X_i|^p] <\infty$ then $E[|\sum_{i=1}^{2^{k+1}N} X_i |^p] \le C (2^{k+1}N)^{p/2}$ for some
$C<\infty$ (see item 16, page 60 of \cite{petrov}), 
and hence for another $C<\infty$,
\begin{eqnarray} &
P[ \sup_{n\ge N}\frac{1}{n}\sum_{i=1}^n X_i \ge \epsilon ] \le   C 
\epsilon^{-p} N^{-p/2}.&\label{kolmogorov}
\end{eqnarray}

Applying (\ref{kolmogorov}) to (\ref{you}), by Lemma \ref{moment1}
we obtain (\ref{such-that}).
 \end{proof}

\subsection{Estimates for the enlarged process.}  We start with a few standard estimates for hitting times of 
random walks.
\begin{lemma}
\label{line-est} Let $\{X_t:t\ge 0\}$ be a continuous time simple
symmetric random walk on ${\mathbb Z}$, of total jump rate $2$, starting
from $x\le -1$, $c>0$, and
\begin{equation}
\label{tau}
\tau:=\inf\{t\ge 0: X_t\ge \lfloor ct\rfloor\}.
\end{equation}
Then \begin{equation} \nonumber
P[\tau=\infty]\ge \begin{cases} 1-\exp\{(1+x)\theta_c\} & x\le -2 \\  \exp\{-2/c\}(1-\exp\{-\theta_c\}) & x=-1\end{cases}
\end{equation}
where $\theta_c>0$ is the nonzero solution of  $c\theta -2(\cosh\theta -1)=0$,
and there exist $0< C,C'<\infty$ such that 
\begin{equation}\label{cccc}
P[t< \tau<\infty] < C'\exp \{ -C
t\}\Big(\exp\{-\theta_c|x|/2\}+\exp\{-tI\big(|x|/(2t)\big)\} 
\Big),
\end{equation}
where $I(u)=2+u\sinh^{-1}(u/2)-\sqrt{4+u^2}$ is the rate function for a random walk.
\end{lemma}

\begin{proof}[Proof] For $\theta\in \mathbb R$,
$\exp\{\theta X_t-2(\cosh\theta -1)t\}$
is a martingale. By the optional stopping time theorem, 
\begin{equation}
\label{it-follows}
E[\exp\{\theta X_{\tau\land n}-2(\cosh\theta -1)\tau\land n\}]=\exp\{\theta x\}.
\end{equation}
Now, $X_{\tau\land n}\le \lfloor c\tau\land n\rfloor$, hence
$\theta X_{\tau\land n}-2(\cosh\theta -1)\tau\land n
\le (c\theta -2(\cosh\theta -1))\tau\land n$.
It follows that if $\theta\ge \theta_c$, we can apply the bounded
convergence theorem in equality (\ref{it-follows}) taking the
limit when $n\to\infty$, to conclude that,
\begin{equation}
\label{conclude-that}
E[{\bf 1}(\tau<\infty)\exp\{\theta X_\tau-2(\cosh\theta -1)\tau\}]=\exp\{\theta x\}.
\end{equation}
When $\tau<\infty$, we have  $X_\tau=\lfloor c\tau\rfloor$ and therefore 
$\theta X_\tau-2(\cosh\theta -1)\tau\ge
(c\theta -2(\cosh\theta -1)) \tau-\theta$.
Letting  $\theta\downarrow\theta_c$ in (\ref{conclude-that})  and
applying the bounded convergence theorem, we obtain,
\begin{equation} \nonumber
P[\tau<\infty]\le \exp\{(1+x)\theta_c\} .
\end{equation}
To get the bound for $x=-1$, note that the probability that the random walk does not move before time $t = 1/c$
 is $\exp\{-2/c\}$.  Then use the Markov property and the result for $x=-2$.

To prove (\ref{cccc}):   \begin{equation}
P[t<\tau<\infty] \le P[ X_t > B] + P[ t<\tau<\infty, X_t \le B]  
\end{equation}
Now
\begin{equation} \nonumber
P[X_t >B] \le \exp \{ - t I( (B+|x|)/t) \}
\end{equation}
and  by the strong Markov property and translation invariance,

\begin{equation} \nonumber
P[ t<\tau<\infty, X_t \le B] \le P_{-(\lfloor ct\rfloor-B)}[ \tau<\infty]
  \le \exp \{ -(\lfloor ct\rfloor-B-1)\theta_c\}. 
\end{equation}
Choosing $B=(\lfloor ct\rfloor +x)/2$, on the above inequalities,
and using the convexity of the rate function $I$, gives (\ref{cccc}). 
\end{proof}

\begin{lemma} 
\label{est-V}  For each $\alpha'<\alpha$, there exists $C<\infty$ such that for $t\ge 1$, and all $w\in\mathbb S$,
\begin{equation}
 {\mathbb P}_w [ t<V<\infty ]\le
C\exp \{-C t \}.
\end{equation}
\end{lemma}
\begin{proof} Note that in the worst case in which there are $M$ random walks at each site
to the left of the origin, we get the bound,

\begin{equation} \nonumber
{\mathbb P}_w[t<V<\infty]\le M\sum_{x=-1}^{-\infty} P_x[t<\tau<\infty],
\end{equation}
where $\tau$ is defined in display (\ref{tau}).
On the other hand, it is true that $\sum_{k=1}^\infty e^{-I(k/(n+1))}\le (n+1)\sum_{k=0}^\infty e^{-I(k)}$.
This estimate, the previous inequality and inequality (\ref{cccc}) of Lemma \ref{line-est} give us the result.
\end{proof}

\medskip
 From Lemmas \ref{est-U} and \ref{est-V} we have
\begin{corollary} 
\label{est-D} For each $0\le p < M/2$ there is a $C<\infty$ depending only on $p$, $M$ and $\alpha'$
 such that for all initial conditions $w$,
\begin{equation} \nonumber
{\mathbb P}_w\left[ t<D<\infty\right]\le Ct^{-p/2}.
\end{equation}
\end{corollary}

\begin{lemma}
\label{delta1}
    There is a $\delta_1>0$ such that,  
\label{est-V-infty}
\begin{equation} \nonumber
{\mathbb P}_w\left[V<\infty\right]<1-\delta_1.
\end{equation}
\end{lemma}
\begin{proof}[Proof] Without loss of generality, $r=0$.
Now, take the worst case scenario where $w$ has
 $M$ particles at each site $x\le -1$. By Lemma \ref{line-est},
\begin{equation} \nonumber
{\mathbb P}_w[V=\infty]\ge e^{-2M/\alpha'}(1-e^{-\theta_{c}})^M
\Pi_{n=1}^\infty (1-e^{-n\theta_{c}})^M=\delta_1>0.
\end{equation}
\end{proof}

\begin{lemma}
\label{delta2}   Suppose that $M>4$. There is a $\delta_2>0$ such that for all initial conditions $w$ with at least $a-1$ particles at the rightmost
visited site $r$  
\label{est-U-infty}
\begin{equation} \nonumber
{\mathbb P}_w\left[U<\infty\right]<1-\delta_2.
\end{equation}
\end{lemma}
 
\begin{proof}[Proof]  We can also assume that $r=0$. To estimate ${\mathbb P}_w[U=\infty]$ below, note that
\begin{eqnarray}
\label{now}&
{\mathbb P}_w[U=\infty]=P[\cap_{k=1}^\infty\{\sum_{j=1}^k\nu_j< k/\alpha'\}].&
\end{eqnarray}
Let $n\in{\mathbb N}$ and $0<\epsilon<1/(2\alpha')$ 
and define $G$ to be the event that
each random walk $Y_{(x,i)}$ with $0\le x\le n$, moves
$M+1$ steps to the right before time $\epsilon$.
When $G$ happens, $\nu_k<1/(2\alpha')$ and hence $\sum_{j=1}^k\nu_j<k/(2\alpha')$
for all
$k\le n+\lfloor M/a\rfloor:=n'$.  Note that $G\cap \cap_{k=1}^\infty\{\sum_{j=1}^k\nu_j< k/\alpha'\} \supset G \cap H$ where
\begin{eqnarray}
\label{now2}&
 H:= \bigcap_{k=n'+1}^\infty
\left\{\sum_{j=n'+1}^k\nu_j<k/\alpha' - n'\epsilon
\right\}&
\end{eqnarray}
Furthermore, $G$ and $H$ are independent so ${\mathbb P}_w[U=\infty]\ge P[G]P[H]$.
Now
\begin{equation} \nonumber
P[H^c] \le \sum_{k=n'+1}^\infty  P[ \sum_{j=n'+1}^k\nu_j\ge k/\alpha' -n'\epsilon]
\end{equation}
and letting $\gamma_j = \nu_j- E[\nu_j] = \nu_j-1/\alpha$, for $2\le p<M/2$, with $\beta :=1/\alpha'- 1/\alpha>0$
\begin{eqnarray} &
 P[ \sum_{j=n'+1}^k\nu_j\ge k/\alpha' -n'\epsilon]\le E[ (\sum_{j=n'+1}^k\gamma_j)^p]  (k\beta -n'\epsilon)^{-p}&\nonumber \\ &
\le C k^{p/2}(k\beta -n'\epsilon)^{-p}, &
\end{eqnarray}
where in the last inequality we used the same estimates explained between
displays (\ref{thirtysix}) and (\ref{kolmogorov}).
Taking $\epsilon= \beta /(2n')$ gives $k^{p/2}(k\beta -n'\epsilon)^{-p} \le C' k^{-p/2}$.  Hence
$
P[H^c] \le C\sum_{k=n'+1}^\infty  k^{-p/2} 
$ so as long as $p>2$ (which is possible since $M>4$) we obtain that $P[H^c] < 1- \delta_2<1$ for sufficiently large $n$.
Choose such an $n<\infty$, and note that for this $n$, $P[G] \ge \delta_3>0$ as well.\end{proof}

\begin{lemma}   Suppose that $M>4$. There is a $\delta>0$ such that for all initial conditions $w$ with at least $a-1$ particles at the rightmost
visited site $r$  
\label{est-D-infty}
\begin{equation}
{\mathbb P}_w\left[D<\infty\right]<1-\delta.
\end{equation}
\end{lemma}

\begin{proof}[Proof]   Since $U$ and $V$ are independent by part $(iii)$ of Lemma \ref{coupling},
${\mathbb P}_w[D<\infty]=1-{\mathbb P}_w[U=\infty]{\mathbb P}_w[V=\infty]$.

\end{proof}

\begin{lemma}\label{speed}  There is a $C<\infty$ such that 
 for every initial condition $w$ with $r=0$, and any $M'>M$,
\begin{equation}
{\mathbb P}_w\left[ r_t\ge M' t\right]\le C\exp\{-Ct\}.
\end{equation}
\end{lemma}

\begin{proof}  Note that $r_t$ is a process on $\mathbb Z$, increasing by one whenever a  particle
jumps to the right from $r$.  Since there at most $M$ particles there, the maximum jump rate is $M$.
The Lemma then follows from standard estimates on Poisson processes.
\end{proof}

\begin{lemma} \label{speed2} For each $p<M/2$ there is a $C<\infty$ such that
\begin{equation}\label{tailk}
{\mathbb P}_{a\delta_0}[\kappa_1>t|U=\infty]\le Ct^{-p/2}.
\end{equation}
\end{lemma}

\begin{proof} Let us first write,
\begin{equation} \nonumber
{\mathbb P}_{a\delta_0}\left[\left. \kappa_1>t \right| U=\infty\right]
=\sum_{k=1}^\infty {\mathbb P}_{a\delta_0}\left[\left. S_k >t, K=k \right| U=\infty\right].
\end{equation}
Applying recursively the strong Markov property
to the stopping times $\{S_j:j\ge 1\}$ we see that for every $k\ge 1$,
\begin{equation} \nonumber
{\mathbb P}_{a\delta_0}\left[\left. S_k >t, K=k \right| U=\infty\right]
\le (1-\delta)^{k-1},
\end{equation}
where $\delta>0$ is given by Lemma
 \ref{est-D-infty}.   For any $\ell>0$ we therefore have,
\begin{eqnarray}
\label{decomp}&
{\mathbb P}_{a\delta_0}\left[\left. \kappa_1>t \right| U=\infty\right]
\le
\sum_{k=1}^{\ell} 
{\mathbb P}_{a\delta_0}\left[\left. t<S_k <\infty \right| U=\infty\right]
+\delta^{-1} (1-\delta)^\ell.&
\end{eqnarray}
Let $1>\gamma>0$ and consider the event
\begin{equation} \nonumber
A_k = \{ r_{D_1}- r_{S_1} < t^\gamma, r_{D_2}- r_{S_2} < t^\gamma, \ldots, r_{D_{k-1}}- r_{S_{k-1}} < t^\gamma\}
\end{equation}
On $A_k$ we have
$r_{S_k} \le k(L +t^\gamma)$.
Since $\tilde r_t \le r_t$, if $U=\infty$, then $r_t >\lfloor \alpha' t \rfloor$ for all $t>0$.  Therefore, on $A_k\cap \{U=\infty\}$,
\begin{equation}
\lfloor \alpha' S_k\rfloor  \le k(L +t^\gamma).
\end{equation}
Hence for $t> (\ell (L+ t^\gamma)+1) / \alpha'$ and $k\le \ell$,
\begin{equation} \nonumber
{\mathbb P}_{a\delta_0}\left[\left. t<S_k <\infty , A_k\right| U=\infty\right]=0
\end{equation}
and therefore 
\begin{eqnarray}&
{\mathbb P}_{a\delta_0}\left[\left. t<S_k <\infty \right| U=\infty\right]
\le {\mathbb P}_{a\delta_0}\left[\left. A_k^c , S_k<\infty\right| U=\infty\right].&\label{lkj}
\end{eqnarray}
By Lemma \ref{est-U-infty}, since ${\mathbb P}_{a\delta_0}[U=\infty]\ge
\delta_2
 >0$.  So for some $C<\infty$, the right hand side of (\ref{lkj}) is bounded above
by
\begin{eqnarray} &
C\sum_{i=1}^{k -1}{\mathbb P}_{a\delta_0}\left[ r_{D_i}- r_{S_i}\ge  t^\gamma, S_k <\infty \right].&
\end{eqnarray}
Now let $M'>M$ and
\begin{eqnarray} &
\!\!\!\!\!{\mathbb P}_{a\delta_0}\!\left[ r_{D_i}- r_{S_i}\ge  t^\gamma, S_k <\infty \right]\! =\! 
{\mathbb P}_{a\delta_0}\!\left[ r_{D_i}- r_{S_i}\ge  t^\gamma, S_k <\infty, D_i-S_i \le t^\gamma/M' \right]
&\nonumber \\&+{\mathbb P}_{a\delta_0}\left[ r_{D_i}- r_{S_i}\ge  t^\gamma, S_k <\infty, D_i-S_i > t^\gamma/M' \right] &\nonumber\\
& \le {\mathbb P}_{a\delta_0}\left[ r_{S_i + t^\gamma/M' }- r_{S_i}\ge  t^\gamma \right]
+{\mathbb P}_{a\delta_0}\left[   t^\gamma/M' <D_i-S_i <\infty\right].&
\end{eqnarray}
Note that in the last equation we used the fact that $S_k<\infty$ implies that $D_i-S_i<\infty$ for $i<k$.
By the strong Markov property and Lemma \ref{speed}, 
\begin{equation} \nonumber
{\mathbb P}_{a\delta_0}\left[ r_{S_i + t^\gamma/M' }- r_{S_i}\ge  t^\gamma \right]\le C\exp{-C t^\gamma}.
\end{equation}
By the strong Markov property and Corollary \ref{est-D}, for each $p<M/2$,
\begin{equation} \nonumber
{\mathbb P}_{a\delta_0}\left[   t^\gamma/M' <D_i-S_i <\infty\right] \le C t^{-\gamma p /2}.
\end{equation}
Choosing $\ell = C\log t$ with $C= p(2\log(1-\delta)^{-1})^{-1}$ from (\ref{decomp}) and the previous estimates
we obtain (\ref{tailk}).
\end{proof}

\begin{corollary}

\label{moment}
 Let $\kappa_1$ be the first regeneration time
of the stochastic combustion process.

\begin{itemize}

\item[a)] For $M> 4$, 
${\mathbb E}_{a\delta_0}\left[\left. \kappa_1\right|U=\infty\right]<\infty, $
and ${\mathbb E}_{a\delta_0}\left[\left. r_{\kappa_1}\right|U=\infty\right]<\infty.$

\item[b)] For $M> 8$,
${\mathbb E}_{a\delta_0}\left[\left. \kappa_1^2\right|U=\infty\right]<\infty,$
and ${\mathbb E}_{a\delta_0}\left[\left. r_{\kappa_1}^2\right|U=\infty\right]<\infty.$

\end{itemize}
\end{corollary}

\begin{proof}[Proof] The statements for $\kappa_1$ follow from 
Lemma \ref{speed2}, and those for $r_{\kappa_1}$  from 
Lemmas \ref{speed} and \ref{speed2}.
\end{proof}

\section{Limit theorems}
\label{limits}
In this section we use the renewal structure to prove Theorem \ref{theorem1}, the law of large
numbers and the central limit theorem for $r_t$.  The argument for the law of large numbers
follows that of Sznitman and Zerner
in \cite{sz}, developed in the context of multi-dimensional transient random
walks in random environments.  The argument for the central limit theorem is from \cite{s1}.

\subsection{ Law of Large Numbers.} We consider the stochastic combustion process started
with an initial condition $r=0$, $\eta(0,0)\in\{1,\ldots, M\}$ and $\eta(0,x)\in\{0,\ldots, M\}$.
We will prove that a.s.
\begin{equation}
\label{as-lim}
\lim_{t\to\infty} \frac{r_t}{t}=v:=\frac{{\mathbb E}_{a\delta_0}[r_{\kappa_1}
|U=\infty]}{{\mathbb E}_{a\delta_0}[\kappa_1|U=\infty]}.
\end{equation}
Let us first note that by Proposition \ref{iid}, we have a.s.,
\begin{equation}
\label{kappa}
\lim_{n\to\infty}\frac{\kappa_n}{n}={\mathbb E}_{a\delta_0}[\kappa_1
|U=\infty],\quad{\rm and}\quad
\lim_{n\to\infty}\frac{r_{\kappa_n}}{n}={\mathbb E}_{a\delta_0}[r_{\kappa_1}
|U=\infty].\end{equation}
For $t\ge 0$, define  $n_t:=\sup\{n\ge 0:\kappa_n\le t\}$, with
the convention that $\kappa_0=0$. Note that (\ref{kappa}) ensures that
$n_t<\infty$ a.s. and by definition we also have,
$
\kappa_{n_t}\le t <\kappa_{n_t+1},
$ and $\lim_{t\to\infty}\kappa_{n_t}=\infty$. It follows from this
inequality and the limit (\ref{kappa}) that, a.s.,
$\lim_{t\to\infty}{n_t}/{t}={1}/{{\mathbb E}_{a\delta_0}[\kappa_1
|U=\infty]}$ and hence almost surely,
\begin{equation}
\label{but-lim}
\lim_{t\to\infty}{r_{\kappa_{n_t}}}/{t}=\lim_{t\to\infty}(r_{\kappa_{n_t}}/{\kappa_{n_t}})( {\kappa_{n_t}}/{t}) =v,
\end{equation}
Now   \begin{equation}
\label{and-lim}
\lim_{t\to\infty}\frac{|r_t-r_{\kappa_{n_t}}|}{t}\le
\lim_{t\to\infty}\frac{r_{\kappa_{n_t+1}}-r_{\kappa_{n_t}}}{t}=0.
\end{equation}
from (\ref{but-lim}).
This proves the law of large numbers.

\subsection{Central limit theorem.}  Starting with the same initial conditions we consider
\begin{equation}
\nonumber
B_t^\epsilon:=\epsilon^{1/2}\left(r_{\epsilon^{-1}{t}}-\epsilon^{-1} {vt}\right),
\qquad\qquad t\ge 0,
\end{equation}
Define
\begin{equation} \nonumber
R_j:=r_{\kappa_{j+1}}-r_{\kappa_j}-(\kappa_{j+1}-\kappa_j)v,\qquad\qquad
j\ge 0,
\end{equation}
and denote for $m\ge 0$ the partial sums
 $\Sigma_m:=\sum_{j=1}^mR_j$. 

For any
$0\le t\le T<\infty$,
\begin{equation} \nonumber
|B_t^\epsilon-\epsilon^{1/2}\Sigma_{n_{t/\epsilon}}
|\le 2\epsilon^{1/2}\sup_{0\le n\le n_{\lfloor \epsilon^{-1}T\rfloor}}(r_{\kappa_{n+1}}-r_{\kappa_n})
+
2v\epsilon^{1/2}\sup_{0\le n\le n_{\lfloor \epsilon^{-1}T\rfloor} }(\kappa_{n+1}-\kappa_n).
\end{equation}
For every $u>0$ we have, by Proposition
\ref{iid}
\begin{eqnarray} &
\nonumber
{\mathbb P}_{a\delta_0}[\sup_{0\le n\le n_{\lfloor \epsilon^{-1}T\rfloor}}
\epsilon^{1/2}(\kappa_{n+1}-\kappa_n)>u]\\
\nonumber
&\le{\mathbb P}_{a\delta_0}[\kappa_1>\epsilon^{-1/2}u]+ u^{-2} \epsilon(
t\epsilon^{-1}+1)
{\mathbb E}_{a\delta_0}[\kappa_1^21(\kappa_1>\epsilon^{-1/2}u)|U=\infty],
\end{eqnarray}
which by part $(b)$ of Corollary \ref{moment}, a.s. converges to $0$ as $\epsilon\to 0$.
Hence, in probability
\begin{equation}
\label{error1}
\sup_{0\le n\le n_{\lfloor \epsilon^{-1}T\rfloor}}
\epsilon^{1/2}(\kappa_{n+1}-\kappa_n)\to 0.
\end{equation}
and similarly
$$
\sup_{0\le n\le n_{\lfloor \epsilon^{-1}T\rfloor}}
\epsilon^{1/2}(r_{\kappa_{n+1}}-r_{\kappa_n})\to 0.$$
Hence, $B_t^\epsilon-\epsilon^{1/2}\Sigma_{n_{\epsilon^{-1}t}}$ converges to
$0$ in probability, uniformly on compact sets of $t$. From Donsker's invariance principle, we know that $\sqrt{\epsilon}\Sigma_{\cdot /\epsilon}$
converges in law to a Brownian motion with variance ${\mathbb E}_{\delta_0}[(r_{\kappa_1}-\kappa_1v)^2|U=\infty]$,
where $\Sigma_s, s\ge 0,$ now stands for the linear interpolation of $\Sigma_m,m\ge 0$.  From the previous proof
we have $\lim_{t\to\infty} n_t/t =1/ {\mathbb E}_{a\delta_0}[\kappa_1|U=\infty]$.  Since
$\epsilon^{-1}k_{t\epsilon^{-1}}$ is increasing in $t$,  the convergence is uniform on compact sets of $t$.  This, together with the convergence in law of $\epsilon^{1/2}\Sigma_{\cdot /\epsilon}$
 to a Brownian motion with variance ${\mathbb E}_{a\delta_0}[(r_{\kappa_1}-\kappa_1 v)^2]$, implies that
 $\epsilon^{1/2}\Sigma_{k_{\lfloor\epsilon^{-1} t\rfloor}}$ is tight in the Skorohod topology, and that its finite dimensional distributions converge to the
finite dimensional distribution of a Brownian motion with variance, 
\begin{equation}
\label{sigma}
\sigma^2:=
\frac{{\mathbb E}_{a\delta_0}[(r_{\kappa_1}-\kappa_1v)^2|U=\infty]}
{{\mathbb E}_{a\delta_0}[\kappa_1|U=\infty]}.
\end{equation}

\subsection{Non-degeneracy of the variance.} It suffices to prove that, for some $\alpha'<\beta<v$,
\begin{equation} \nonumber
{\mathbb P}_{a\delta_0}[r_{\kappa_1}=L,L\beta^{-1}\le \kappa_1~|~ U=\infty]>0.
\end{equation}
Now,
\begin{equation} \nonumber
{\mathbb P}_{a\delta_0}[r_{\kappa_1}=L,L\beta^{-1}\le \kappa_1, U=\infty]\ge
{\mathbb P}_{a\delta_0}[L\beta^{-1}<S_1< U, D\circ \theta_{S_1}=\infty]
\end{equation}
The right hand side we can write as
\begin{equation} \nonumber
{\mathbb E}_{a\delta_0}[1(L\beta^{-1}<S_1<U ){\mathbb E}_{a\delta_0}[  1(V\circ \theta_{S_1}=\infty)1( U\circ \theta_{S_1}=\infty) ~|~{\mathcal F}_{S_1}]]
\end{equation}
Given ${\mathcal F}_{S_1}$, $ V\circ \theta_{S_1}$ and $U\circ \theta_{S_1}$ are independent, so
\begin{eqnarray} &
\mathbb E_{a\delta_0}[  1(V\circ \theta_{S_1}=\infty)1( U\circ \theta_{S_1}=\infty) ~|~\mathcal F_{S_1}]] &\nonumber\\ &
=\mathbb P_{a\delta_0}[  V\circ \theta_{S_1}=\infty ~|~\mathcal F_{S_1}]]\mathbb P_{a\delta_0}[   U\circ \theta_{S_1}=\infty ~|~\mathcal F_{S_1}]].&
\end{eqnarray}
But since by Lemmas \ref{delta1} and \ref{delta2} we have $\mathbb 
P_{a\delta_0}[   U\circ \theta_{S_1}=\infty ~|~\mathcal F_{S_1}]]
=\mathbb P_{a\delta_0}[U=\infty]\ge\delta_2>0$
and $\mathbb 
P_{a\delta_0}[   V\circ \theta_{S_1}=\infty ~|~\mathcal F_{S_1}]]
\ge\delta_1>0$
we have,

$$
{\mathbb P}_{a\delta_0}[L\beta^{-1}<S_1< U, D\circ \theta_{S_1}=\infty]\\
\ge \delta_1\delta_2{\mathbb P}_{a\delta_0}[L\beta^{-1}<S_1< U].
$$
But it is easy to check that 
${\mathbb P}_{a\delta_0}[L\beta^{-1}<S_1< U]>0$.
In fact, it is enough to lower bound this probability by the probability
that one of the random walks at site $0$ moves to site $L$ in
a time $t$ such that $L\beta^{-1}<t<L(\alpha')^{-1}$,
and then stays at site $L$ between time $t$ and time $L(\alpha')^{-1}$,
 while all other random walks between
sites $0$ and $L$ do not move at all during the time
interval $[0,L(\alpha')^{-1}]$.

\end{document}